\documentclass[10pt,english]{amsart}
\usepackage{amsmath}
\usepackage{amsxtra}
\usepackage{amscd}
\usepackage{amsthm}
\usepackage{amsfonts}
\usepackage{amssymb}
\usepackage[matrix,arrow,curve]{xy}
\usepackage{hyperref}
\usepackage{rtr}

\newtheorem{lem}[subsection]{Lemma}

\newtheorem{thm}[subsection]{Theorem}
\theoremstyle{definition}
\newtheorem{defn}[subsection]{Definition}
\theoremstyle{remark}
\newtheorem{rem}[subsection]{Remark}

\let\ncsave\newcommand
\let\newcommand\providecommand

\newcommand{\nc}{\newcommand}
\nc{\renc}{\renewcommand} \nc{\ssec}{\subsection}
\nc{\sssec}{\subsubsection} \nc{\on}{\operatorname}

\nc\ol{\overline} \nc\ul{\underline} \nc\wt{\widetilde}
\nc\tboxtimes{\wt{\boxtimes}} \nc{\alp}{\alpha}

\nc{\ZZ}{{\mathbb Z}} \nc{\NN}{{\mathbb N}} \nc{\CC}{{\mathbb C}}
\nc{\OO}{{\mathbb O}} \renc{\SS}{{\mathbb S}} \nc{\DD}{{\mathbb
D}}

\nc{\Fq}{{\mathbb F}_q} \nc{\Fqb}{\ol{{\mathbb F}_q}}
\nc{\Ql}{\ol{{\mathbb Q}_\ell}} \nc{\id}{\text{id}} \nc\X{\mathcal
X}

\nc{\Hom}{\on{Hom}} \nc{\Lie}{\on{Lie}} \nc{\Loc}{\on{Loc}}
\nc{\Pic}{\on{Pic}} \nc{\Bun}{\on{Bun}} \nc{\IC}{\on{IC}}
\nc{\Aut}{\on{Aut}} \nc{\rk}{\on{rk}} \nc{\Sh}{\on{Sh}}
\nc{\Perv}{\on{Perv}} \nc{\pos}{{\on{pos}}} \nc{\Conv}{\on{Conv}}
\nc{\Sph}{\on{Sph}} \nc{\Sym}{\on{Sym}}
\nc{\BunBb}{\overline{\Bun}_B} \nc{\Buno}{\overset{o}{\Bun}}
\nc{\BunPb}{{\overline{\Bun}_P}}
\nc{\BunBM}{\overline{\Bun}_{B(M)}}
\nc{\BunPbw}{{\widetilde{\Bun}_P}}
\nc{\BunBP}{\widetilde{\Bun}_{B,P}} \nc{\GUb}{\overline{G/U}}
\nc{\GUPb}{\overline{G/U(P)}}
\nc{\iso}{{\stackrel{\sim}{\longrightarrow}}}

\nc{\Hhom}{\underline{\on{Hom}}} \nc\syminfty{\on{Sym}^{\infty}}
\nc\lal{\ol{\lambda}} \nc\xl{\ol{x}} \nc\thl{\ol{\theta}}
\nc\nul{\ol{\nu}} \nc\mul{\ol{\mu}} \nc\Sum\Sigma
\nc{\oX}{\overset{o}{X}{}}

\nc{\M}{{\mathcal M}} \nc{\N}{{\mathcal N}} \nc{\F}{{\mathcal F}}
\nc{\D}{{\mathcal D}} \nc{\Q}{{\mathcal Q}} \nc{\Y}{{\mathcal Y}}
\nc{\G}{{\mathcal G}} \nc{\E}{{\mathcal E}} \nc{\CalC}{{\mathcal
C}}
\nc\Dh{\widehat{\D}}

\nc{\C}{{\mathcal C}} \nc{\K}{{\mathcal K}}
\renewcommand{\H}{{\mathcal H}}

\nc{\T}{{\mathcal T}} \nc{\V}{{\mathcal V}} \renc{\P}{{\mathcal
P}} \nc{\A}{{\mathcal A}} \nc{\B}{{\mathcal B}} \nc{\U}{{\mathcal
U}}
\renewcommand{\L}{{\mathcal L}}
\nc{\Gr}{\on{Gr}}

\nc{\frn}{{\check{\mathfrak u}(P)}}
\nc\f{{\mathfrak f}}
\renewcommand\k{{\mathfrak k}}
\nc{\q}{{\mathfrak q}} \nc{\p}{{\mathfrak p}} \nc{\s}{{\mathfrak
s}} \nc\w{\text{w}}

\nc\mathi\iota \nc\Spec{\on{Spec}} \nc\Mod{\on{Mod}}
\nc{\tw}{\widetilde{\mathfrak t}} \nc{\pw}{\widetilde{\mathfrak
p}} \nc{\qw}{\widetilde{\mathfrak q}} \nc{\jw}{\widetilde j}

\nc{\grb}{\overline{\Gr}} \nc{\I}{\mathcal I}

\nc{\lambdach}{{\check\lambda}} \nc{\Lambdach}{{\check\Lambda}{}}
\nc{\much}{{\check\mu}} \nc{\omegach}{{\check\omega}}
\nc{\nuch}{{\check\nu}} \nc{\etach}{{\check\eta}}
\nc{\alphach}{{\check\alpha}} \nc{\betach}{{\check\beta}}
\nc{\rhoch}{{\check\rho}} \nc{\ch}{{\check h}}

\nc{\Hb}{\overline{\H}}

\usepackage{fnklberg}     

\nc{\MLT}{{\mathcal{M}^\theta(v,w)}}
\nc{\ML}{{\mathcal{M}(v,w)}}
\nc{\ra}{\rightarrow}

\newcommand*{\mb}[1]{\mathbb{#1}}
\newcommand*{\mr}{\mathrm}
\newcommand*{\mc}{\mathcal}

\let\newcommand\ncsave

\def\k{\Bbbk}
\DefOps{\diag \Tr \ad \Br}

\begin{document}

\title[On some asymptotic formulas for curves in positive characteristic. ]%
    {On some asymptotic formulas for curves in positive characteristic}

\author{Dmitry Kubrak}
\address[Dmitry Kubrak]{
State University Higher School of Economics\\
Department of Mathematics\\
{20~Myasnitskaya~St},
Moscow 101000, Russia}
\email{dmkubrak@gmail.com}

\maketitle

\begin{abstract}
We prove some general asymptotic formula for the values of $L$-function of a sequence of constructible $\mb Q_l$-sheaves on curves over $\mb F_q$ with some good asymptotic properties. We also give the asymptotic formula for the number of points on the stack $\mr{Bun}_G$ for asymptotically exact sequence of curves and a split reductive group $G$. In the case of $\mr{GL}_n$ we prove that the same formula holds if we take only semistable points in the account.
\end{abstract}


\section{Introduction}
The asymptotic formulas for various arithmetic invariants have always been of great interest to number theorists, starting from classical Brauer-Siegel theorem. Its natural generalisation is given by Tsfasman-Vl\v{a}duts formula for the asymptotic behavior of class number of a curve. We give generalisations of this formula in two directions: the asymptotic formula for the quasi-residues of a family of $l$-adic sheaves having some good properties and asymptotic formula for the number of points on the stack of $G$-bundles on curve for reductive split G. For $\mr{GL}_n$ we are able to prove analogous formula for the number of semistable bundles.

\section{General asymptotic formula for $l$-adic sheaves}

Here we will desribe the problem we are going to study and prove a general asymptotical formula for what we call \textit{asymptotically exact system} with some additional global assumptions. By a system we mean nothing more than a collection of pairs $(X_i,\mathcal F_i)$, each pair consisting of a curve $X_i$ and a constructible sheaf $\mc F_i$ on it. We are interested in the behavior of quasiresidue of $L$-function of $\mc F_i$ as $i$ tends to infinity. The answer we are looking for should extract the result through some local invariants of the system, namely the Tsfasman-Vl\v{a}duts invariants \cite{TV}. The property of a system being asymptotically exact is exactly the property of these invariants to exist. 

\subsection{Asymptotical exactness} Let a \textit{system} ($\Xi$, $\Phi$) be a sequence of pairs $(X_i,\mc F_i)$, numbered by natural numbers, where $X_i$ is some complete smooth curve over $\mb F_q$ and $\mc F_i\in D^b_c(X_i,\mb Q_l)$ is an element of derived category of constructible $l$-adic sheaves on $X_i$. In all our applications all $\mc F_i$'s will be actual sheaves, but formulas can be proved for complexes without any loss of generality. For simplicity we will still call $\mc F_i$ a sheaf even if it is actually an element of derived category. 

Let $i:x\hookrightarrow X, x=\mr{Spec}(\mb F_{q^n})$ be a closed point of $X$ and let $\mc F$ be a sheaf on $X$. Associated to this point we have a local $L$-factor 

$$
L_x(\mc F, s)=\frac{\prod_{i=2k} \det(1-\mr{Fr}_x q^{-ns}| H^i(i^*\mc F))}{\prod_{i=2k+1} \det(1-\mr{Fr}_x q^{-ns}| H^i(i^*\mc F))}.
$$ 

That is actually the only local data of $\mc F_i$'s we are interested in.

\begin{defn}
\textit{Tsfasman-Vl\v{a}duts (TV) structure }$\Upsilon$ on a system ($\Xi$, $\Phi$) is a collection of subsets $Z_{r,i}=Z_r(X_i)\subset|X_i|$, such that
\begin{itemize}
\item every $Z_{r,i}$ is finite,
\item $|X_i|=\sqcup_r Z_{r,i}$,
\item all the points in $Z_{r,i}$ are of the same degree,
\item if $x_1\in Z_{r,i}$ and $x_2\in Z_{r,j}$ , then there is an isomorphism $(\mc F_i)_{x_1}=(\mc F_j)_{x_2}$ as Frobenius-modules.
\end{itemize}
\end{defn}
Morally TV-structure is just a particular choice of grouping points where local $L$-factors of corresponding sheaves are the same. Here is the most natural example of a system with a TV-structure:

\medskip
\noindent\textit{Example 1:} Let $X_i$ be arbitrary and $\mc F_i=\ul{\mb Q_l}$ be a constant sheaf on $X_i$. Then grouping the points by degrees over $\mb F_q$  defines a TV-structure $Z_{r,i}=\{x\in |X|,\deg x =r \}$. Further we will call this system $(\Xi, \ul{\mb Q_l})$.

\medskip
\noindent\textit{Example 2:} Let $X$ be a curve and $\mc F$ be a local system on some open part $j:U\hookrightarrow X$. Let $S$ denote the complement $X\setminus U$. Choose any local system $\mc L$ on $U$ and a sequence of ramified coverings $\varphi_i: X_i\rightarrow X$. By such data we can construct system $(\Xi, \Phi)=\{(X_i, \varphi_i^* j_{!}\mc L)\}$. Here we can take also the full $\mb R^\bullet j_{!}$ or any of its sheaf cohomologies $\mb R^n j_{!}$; by our convention $j_{!}$ applied to an actual sheaf will usually mean the $\mb R^0 j_{!}$. Then by base change the stalks of $\mc F_i$ are easy to describe: for a point $x\in |X_i|$ they depend only on the the image $\varphi_i(x)$ and $\deg x$. So, taking by index set $|X|\times \mb N$, we can define a TV-structure in the following way: $Z_{y,r,i}=\{x\in |X|, \varphi_i(x) = y, \deg x =r \}$.

The next definition is made to distinguish systems with good asymptotic properties: 

\begin{defn}
A system $(\Xi, \Phi)$ with TV-structure $\Upsilon$ is called \textit{asymptotically exact}, if 
\begin{itemize} 
\item the genus $g_{X_i}\ra\infty$ as $i\ra\infty $
\item for every $r$ there exists limit 

$$
\gamma_r(\Xi, \Phi)=\lim_{i\ra\infty} \frac{|Z_{r,i}|}{g_{X_i}}
$$

\end{itemize}

\end{defn}

We call these limits \textit{Tsfasman-Vl\v{a}duts invariants} of a system. They are direct generalisations of classical Tsfasman-Vl\v{a}duts invariants  in the case of  $(\Xi, \ul{\mb Q_l})$, $\gamma_r$ are exactly $\beta_r$ in \cite{TV}. For index $r$, by ${L(r)}$ we will denote the value $L_x(F_i,1)$ for any $x$ in group of points indexed by $r$.

\subsection{Quasi-residue of $L$-function}

Let $\mc F$ be a constructible $l$-adic sheaf on complete curve $X$. $L$-function $L(\mc F, s)$ is defined as a product 
$$
L(\mc F, s)=\prod_{x\in |X|} L_x(\mc F, s)^{-1}=\prod_{x\in |X|} \frac{\prod_{i=2k+1} \det(1-\mr{Fr}_x q^{-s \deg x}| H^i(i_x^*\mc F))}{\prod_{i=2k} \det(1-\mr{Fr}_x q^{-s\deg x}| H^i(i_x^*\mc F))},
$$

where $i_x:x\hookrightarrow X$ is the embedding of the point $x$. By Weil conjectures $L$-function is a rational function of $q^{-s}$:
$$
L(\mc F, s)=\frac{\prod_{i=2k+1} \det(1-\mr{Fr}\cdot q^{-s}| H^i(X,\mc F))}{\prod_{i=2k} \det(1-\mr{Fr} \cdot q^{-s}| H^i(X,\mc F))}.
$$

 In particular there is some integer $r$, such that there exist nonzero limit 
$$
\rho_{(X,\mc F)}=\lim_{s\ra 1}\frac{\L(\mc F,s)}{(s-1)^r},
$$
which we call \textit{quasi-residue} of $\mc F$.

For a constant sheaf we have $r=1$ and

$$
\rho_{(X,\ul{\mb Q_l})}= \frac{q^{(1-g(X))} \cdot \#|\Pic^0(X)|}{q-1}
$$

\subsection{Asymptotic formula for $\rho_{\mc F}$}

Let $(\Xi, \Phi)$ be a system. We are interested in the asymptotical behavior of the quasi-residue as $i$ grows, namely we want to find limit
$$
\lim_{i\ra\infty}\frac{\log_q \rho_{\mc F_i}}{g_{X_i}}.
$$

In the case of $(\Xi, \ul{\mb Q_l})$ which is asymptotically exact, the answer can be deduced from Tsfasman-Vl\v{a}duts formula \cite{TV}
$$
\lim_{i\rightarrow \infty } \frac{\log_{q} \rho_{X_i,\ul{\mb Q_l}}}{g_{X_i}} = -\sum_{m=1}^{\infty} \beta_{m}(\tilde{\mc X}) \log_{q} \left(\frac {q^{m}-1}{q^{m}}\right).
$$

Now we are going to prove a general formula for case of asymptotically exact system satisfying some good properties, listed below

\begin{defn}
\label{localass}
We say that $(\Xi, \Phi)$ satisfy \textit{local assumptions} if the following properties are satisfied

\begin{itemize} 
\item there exist natural number $n$, such that for any $r$ and $i>n$ all Frobenius-weights of stalks $ (\mc F_i)_x$ are less or equal than  $\frac 1 2$. 
\item the total dimension  of the stalk $d((\mc F_i)_x)=\sum \dim H^i((\mc F_i)_x))$ is uniformly bounded by some number $d$ (meaning $d((\mc F_i)_x)<d$ for any $x\in X_i$).

\end{itemize}
\end{defn}

\begin{rem}
As we are interested only in asymptotic behavior,  throwing out all $(X_i,\mc F_i)$ with $i<n$ we get the system with the same asymptotic properties, so without loss of generality we can assume that condition from \ref{localass} is satisfied for all $i$
\end{rem}

For $H^\bullet(X_i,\mc F_i)$ let $H^\bullet_{< 1}(X_i,\mc F_i)$ denote the part of weight less than 1 with respect to Frobenius operator. Let also  $H^\bullet_{\ge 1}(X_i,\mc F_i)$ be the transversal part of weight $\ge 1$.

\begin{defn}
\label{globalass}
We say that $(\Xi, \Phi)$ satisfy \textit{global assumptions} if the following properties are satisfied
\begin{itemize}
\item $H^\bullet(X_i,\mc F_i)$ is a uniformly bounded graded space, namely there exist $n$, such that for any $i$ $H^m(X_i,\mc F_i)>0$, if $|m|>n$,
\item  the weights of $H^\bullet(X_i,\mc F_i)$ are uniformly bounded from above,
\item the total dimension of $H^\bullet_{ \ge 1}(X_i,\mc F_i)$ is ${O}(1)$,\item the total dimension of $H^\bullet_{ <1}(X_i,\mc F_i)$ is ${O}(g_{X_i})$.
\item there exist a weight $\omega<1$, such that for any $i$ all weights of $H^\bullet_{ <1}(X_i,\mc F_i)$ are less than $\omega$.

\end{itemize}
\end{defn}

\begin{thm} (General asymptotic formula).
\label{general} Let $(\Xi, \Phi)$ be an asyptotically exact system with TV-structure $\Upsilon$, satisfying local and global assumptions. Then we have the following formula: 
$$
\lim_{i\rightarrow \infty } \frac{\log_{q} \rho_{\mc F_i}}{g_{X_i}} = -\sum_r \gamma_{r}(\Xi, \Phi) \log_{q}{L(r)}.
$$

\begin{rem}
for $(\Xi, \ul{\mb Q_l})$ we get classical Tsfasman-Vl\v{a}duts formula.
\end{rem}

\end{thm}

\textit{Proof:}
 The proof follows essentially the original proof of Tsfasman and Vl\v{a}duts, but being adobted to the case of any sheaf.  

\begin{lem}
If $(\Xi, \Phi)$ with TV-structure $\Upsilon$ is asymptotically exact, then so does $(\Xi, \ul{\mb Q_l})$
\end{lem}

\textit{Proof:}
Just group $Z_{r,i}$ by the degree of points: let $B_m(X_i)$ be the set of points of degree $m$, then $B_m(X_i)=\cup_{r\in R'} Z_{r,i}$ for some $R'$ and 
$$
\beta_m(\Xi)=\lim_{i\rightarrow \infty } \frac{|B_m(X_i)|}{g_{X_i}}=\sum_{r\in R'} \gamma(\Xi, \Phi).
$$

The key point of the proof is Lefshetz fixed point formula: $\sum_i (-1)^i \mr{Tr}(\Fr^m q^{-m}| H^i(X,\mc F))= \sum_{x\in X(\mb F_{q^m})} \mr{Tr}(\Fr^m q^{-m}|\mc F_x)$. If a point $x\in X(\mb F_{q^m})$ actually came from $x'$ defined over $\mb F_{q^n}$, then $\mr{Tr}(\Fr^m q^{-m}|\mc F_x)=\mr{Tr}(\Fr_{x'}^{\frac{m}{n}}q^{-m}|\mc F_x')$ and there would be exactly $n$ such summands (for all points got from $x'=\mr{Spec} \ \mb F_{q^n}$ by base change to $\mb F_{q^m}$). So we can regroup the summands in Lefshetz formula in the following way:
$$
\sum_j (-1)^j \mr{Tr}(\Fr^m q^{-m}| H^j(X_i,\mc F))=\sum_{n|m} n \sum_{Z_{r,i}\subset B_n}|Z_{r,i}|\ \mr{Tr}(\Fr_{x}^{\frac{m}{n}}q^{-m}|\mc F_x,\ x\in Z_{r,i})
$$.

For simplicity we will denote $\mr{Tr}(\Fr_{x}^{\frac{m}{n}}q^{-m}|\mc F_x,\ x\in Z_{r,j})$ by $\mr{Tr}_{r,m}$.

Let $f(g): \mb N\ra\mb N$ be a function such that $f(g)\ra \infty$ as $g\ra \infty$ , and 

$$
\lim_{g\ra\infty}\frac{f(g)}{\log g}=0.
$$

Let it in addition grow slow enough to satisfy the following property:
for any $m<f(g)$ and any $r$, such that $Z_{r,i}\subset B_m(X_i)$, the differnce 

$$
|\gamma_{r}(\Xi, \Phi)-\frac{|Z_r(X_i)|}{g_{X_i}}|\le \frac{\epsilon(g_{X_i})}{|B_{\le m}(X_i)|}
$$

for some positive real-valued function $\epsilon(g)\ra 0$ as $g\ra \infty$. It is easy to see that such function $f$ can be found easily for any countable set of converging sequences.

\begin{lem}
\label{lefshetz}
For such function $f(g)$ 

$$
\sum_{m=1}^{f(g_{X_i})}\frac{1}{m}\sum_{n|m} n \sum_{Z_{r,i}\subset B_{\le n}} \mr{Tr}_{r,m} \cdot|Z_r(X_i)| = \sum_{m=1}^{f(g_{X_i})} \sum_{Z_{r,i}\subset B_{ m}}\log_q(L(r))\cdot|Z_r(X_i)|+o(g_{X_i})
$$

\end{lem}
\textit{Proof:}
Let's change the summation order from the left side:
$$
\sum_{m=1}^{f(g_{X_i})}\frac{1}{m}\sum_{n|m} n \sum_{Z_{r,i}\subset B_{\le n}} \mr{Tr}_{r,m} \cdot|Z_r(X_i)| = \sum_{n=1}^{f(g_{X_i})}\sum_{Z_{r,i}\subset B_{\le n}}|Z_r(X_i)|\sum_{m=1}^{\lfloor\frac{f(g)}{n}\rfloor+1} \frac{1}{m}\mr{Tr}_{r,mn}. 
$$
Let's also expand the logarithm on the right side (which makes sense because of \ref{localass}):
$$
\sum_{m=1}^{f(g_{X_i})} \sum_{Z_{r,i}\subset B_{\le n}}\log_q(L(r))\cdot|Z_r(X_i)|=\sum_{m=1}^{f(g_{X_i})}\sum_{Z_{r,i}\subset B_{ m}}|Z_r(X_i)|\sum_{n=1}^\infty \frac{1}{n}\mr{Tr}_{r, nm}
$$
We see that the difference is equal to 
$$\sum_{m=1}^{f(g_{X_i})}\sum_{Z_{r,i}\subset B_{ m}}|Z_r(X_i)|\sum_{n=\lfloor\frac{f(g)}{m}\rfloor}^\infty \frac{1}{m}\mr{Tr}_{r,mn}.
$$
Now, for any $r$, such that $Z_{r,i}\subset B_{\le n}$, we have
$$
|\Tr_{r,nm}|=|\mr{Tr}(\Fr_{x}^{\frac{m}{n}}q^{-nm}|\mc F_x,\ x\in Z_{r,\bullet})|\le d(\mc F_x)q^{-nm}< dq^{-\frac{nm}{2}}
$$
using \ref{localass}. So, 
\begin{eqnarray*}
\Bigg|\sum_{m=1}^{f(g_{X_i})}\sum_{Z_{r,i}\subset B_{ m}}|Z_r(X_i)|\sum_{n=\lfloor\frac{f(g)}{m}\rfloor+1}^\infty \frac{1}{m}\mr{Tr}_{r,mn}\Bigg|&\le& d\sum_{m=1}^{f(g_{X_i})} \sum_{Z_{r,i}\subset B_{ m}}|Z_r(X_i)|\sum_{n=\lfloor\frac{f(g)}{m}\rfloor+1}^\infty q^{-\frac{nm}{2}}\le \\ &\le& 3d\sum_{m=1}^{f(g_{X_i})}\sum_{Z_{r,i}\subset B_{ m}}|Z_r(X_i)|q^{-\lfloor\frac{f(g)}{m}\rfloor m-m} =\\ &=& 3d \sum_{m=1}^{f(g_{X_i})} |B_m(X_i)|q^{-\frac{1}{2}\lfloor\frac{f(g)}{m}\rfloor m-\frac{m}{2}} <\\
&<& 3d \sum_{m=1}^{f(g_{X_i})} |B_m(X_i)|q^{-\frac{m}{2}-\frac{f(g)}{4}}
\end{eqnarray*}
as $\lfloor\frac{f(g)}{m}\rfloor m\ge \frac{f(g)}{2}$ for $m\le f(g)$.

In \cite{TV} Tsfasman and Vl\v{a}duts have proven that for any $f(g)$, such that  $\frac{f(g)}{\log g}\ra 0$

$$
\overline{\lim_{i\rightarrow \infty}} \,\frac{1}{g_{X_i}}\,\sum_{m=1}^{f(g_{X_i})}\frac{mB_{m}(X_i)}{q^{m/2 }-1}\le 1,
$$
so
$$
 \lim_{i\rightarrow\infty}q^{-\frac{f(g)}{4}}\frac{3d}{g_{X_i}}\sum_{m=1}^{f(g_{X_i})}B_m(X_i)q^{-\frac{m}{2}}\le \lim_{i\rightarrow\infty} \frac{3d}{g_{X_i}}=0.
$$

Now we will prove that left and right sides of equation in \ref{lefshetz} asymptotically coincide with the left and right sides of \ref{general}.

\begin{lem}
\label{left}
$$
\lim_{i\rightarrow \infty } \frac{\log_{q} \rho_{\mc F_i}}{g_{X_i}} = \lim_{i\rightarrow \infty }\frac{1}{g_{X_i}}\sum_{m=1}^{f(g_{X_i})}\frac{1}{m}\sum_j (-1)^j \mr{Tr}(\Fr^m q^{-m}| H^j(X_i,\mc F_i)).
$$

\end{lem}
\textit{Proof:}
At first let's decompose $\rho_{\mc F_i}$ as a product $\rho_{\mc F_i}^+ \rho_{\mc F_i}^-$, where multiplicants correspond to $H^\bullet_{ \ge 1}(X_i,\mc F_i)$ and $H^\bullet_{ < 1}(X_i,\mc F_i)$. Then, as the total dimension of $H^\bullet_{ \ge 1}(X_i,\mc F_i)$ is bounded and so are the weights, we conclude that $\lim_{i\rightarrow \infty } \frac{\log_{q} \rho_{\mc F_i}^+}{g_{X_i}}=0$. So $\lim_{i\rightarrow \infty } \frac{\log_{q} \rho_{\mc F_i}}{g_{X_i}}=\lim_{i\rightarrow \infty } \frac{\log_{q} \rho_{\mc F_i}^-}{g_{X_i}}$. In $\rho_{\mc F_i}^-$ all weights are $< 1$, so we can expand the logarithm:  

$$
\log_q\rho_{\mc F_i}^-=\sum_{m=1}^{\infty}\frac{1}{m}\sum_j (-1)^j \mr{Tr}(\Fr^m q^{-m}| H^j_{\le 1}(X_i,\mc F_i)).
$$
Now it is left to prove two facts: that 
$$
\lim_{i\rightarrow \infty } \frac{1}{g_{X_i}}\sum_{m=f(g)+1}^{\infty}\frac{1}{m}\sum_j (-1)^j \mr{Tr}(\Fr^m q^{-m}| H^j_{< 1}(X_i,\mc F_i))=0
$$
and that 
$$
\lim_{i\rightarrow \infty } \frac{1}{g_{X_i}}\sum_{m=1}^{f(g)}\frac{1}{m}\sum_j (-1)^j \mr{Tr}(\Fr^m q^{-m}| H^j_{\ge 1}(X_i,\mc F_i))=0
$$.

Expression in the first one can be bounded by geometric progression and the total dimension of weight$<1$ part: 
\begin{eqnarray*}
\lim_{i\rightarrow \infty } \frac{1}{g_{X_i}}\sum_{f(g)+1}^{\infty}\frac{1}{m}\sum_j (-1)^j \mr{Tr}(\Fr^m q^{-m}| H^j_{< 1}(X_i,\mc F_i))\le \lim_{i\rightarrow \infty } \frac{d(H^\bullet_{< 1}(X_i,\mc F_i))}{g_{X_i}} \frac{q^{-(1-\omega)f(g_{X_i})}}{q^{(1-\omega)}-1}=0
\end{eqnarray*}
as $f(g_{X_i})\ra \infty$ and $d(H^\bullet_{< 1}(X_i,\mc F_i))$ is $O(g_{X_i})$. 

The second follows just from the uniformal boundness of weights and the fact that $d(H^\bullet_{\ge 1}(X_i,\mc F_i))$ is $ O(1)$: let $\alpha$ be the bound for weights, then
\begin{eqnarray*}
\lim_{i\rightarrow \infty } \frac{1}{g_{X_i}}\sum_{m=1}^{f(g)}\frac{1}{m}\sum_j (-1)^j \mr{Tr}(\Fr^m q^{-m}| H^j_{\ge 1}(X_i,\mc F_i))&\le& \lim_{i\rightarrow \infty }\frac{d(H^\bullet_{\ge 1}(X_i,\mc F_i))}{g_{X_i}}\sum_{m=1}^{f(g)}\frac{1}{m}q^{(\alpha -1)m}\le\\&\le&\lim_{i\rightarrow \infty }d(H^\bullet_{\ge 1}(X_i,\mc F_i))\frac{q^{(\alpha -1)(f(g)+1)}}{g_{X_i}}=\\
&=&0,
\end{eqnarray*}
as $\lim_{g\ra\infty}\frac{f(g)}{\log g}=0$. The statement of the lemma follows.

The final step of the proof is given by this last lemma
\begin{lem}
\label{right}
$$
\lim_{i\rightarrow \infty } \frac{1}{g_{X_i}}\sum_{m=1}^{f(g_{X_i})} \sum_{Z_{r,i}\subset B_{ m}}\log_q(L(r))\cdot|Z_r(X_i)|=\sum_r \gamma_{r}(\Xi, \Phi) \log_{q}{L(r)}.
$$

\end{lem}

\textit{Proof:}
Again due to \ref{localass} we know that for $Z_{r}(X_i)\subset B_m(X_i)$
$$
\log_{q}{L(r)}\le d\log_q(1-q^{-\frac{m}{2}})\le 3dq^{-\frac{m}{2}}
$$
Remembering the properties of $f(g)$ we get 
\begin{eqnarray*}
&\Big| & \frac{1}{g_{X_i}}\sum_{m=1}^{f(g_{X_i})} \sum_{Z_{r,i}\subset B_{ m}}\log_q(L(r))\cdot|Z_r(X_i)|-\sum_r \gamma_{r}(\Xi, \Phi) \log_{q}{L(r)}\Big|\le\\ &&\le  3\epsilon(g_{X_i})d\sum_{m=1}^{f(g_{X_i})}q^{-\frac{m}{2}}\ra 0, \mbox{as } i\ra\infty
\end{eqnarray*}
This ends the proof of the lemma.

It's easy to see that three proven lemmas (\ref{lefshetz},\ref{left},\ref{right}) give the statement of the theorem.

\section{Asymptotic formulas for $G$-bundles}
Here we will apply the results of previous section to find the asymptotic behaviour of the number of points on stacks $\mr{Bun}_G$ for $G$ split. Using Zagier's formula for the number of semistable bundles (for $\mr{GL}_n$) we also prove that the answer does not change if we restrict ourselves only to semistable part of $\mr{Bun}_{\mr{GL}}$.
\subsection{Asymptotic formula for $|\mr{Bun}^0_G(\mb F_q)|$}
For split reductive group $G$, $\mr{Bun}_G$ as usual denotes the stack of $G$-bundles on  curve $X$, $\mr{Bun}^0_{G}$ denotes the connected component of the  trivial bundle. By $|\mr{Bun}^0_G(\mb F_q)|$ we mean the number of its points in the stacky sense, namely 

$$
|\mr{Bun}^0_G(\mb F_q)|= \sum_{x\in \mr{Bun}^0_G(\mb F_q)}\frac{1}{|\mr{Aut}x|}.
$$
We call it the $G$-mass of $X$ and further denote by $M_G(X)$.

To introduce the formula for $M_G(X)$, we need at first to express it through the quasi-residues of some $L$-functions. 

Let $G$ be a reductive group over $\mb F_q$ of dimension d, with Borel subgroup  $B\subset G$, defined over $\mb F_q$, with maximal torus $T\subset B$ also defined over $\mb F_q$. Let $V=S^{\bullet}(X(T)\otimes \mb Q)$ and $W=(N_G(T)/T)(\ol {\mb F}_q)$. Then $V=S^{\bullet}(X(T)\otimes \mb Q)^W$ is an algebra of polinomials $S^{\bullet }(V_1\oplus \ldots \oplus V_{\mr {dim}  G/B})$ and each $V_i$ is a module over $\mr{Gal}(\ol {\mb F}_q/\mb F_q)$. Then we have Steinberg's formula:
 
$$
|G(\mb F_q)| =q^d\prod_{i\geq1} \det(1-q^{-i}\mr{Fr}|V_i).
$$
 
Let also $\widehat{G}$ be the group of characters $\widehat G\hookrightarrow X(T)$. Then we have $V_1 \cong \widehat G$. Each of $V_i$ tautologically defines a sheaf on $\Spec \mb F_q$, and by pull-back also a sheaf $\mc V_i = \mr{pr}_X^* V_i$ on $X$.

Let $G$ be a reductive group over $\mb F_q(X)$. Choose any left invariant differential top form $\omega$ on $G$. Then $\omega$ defines top-forms $\omega_x$ on $G(F_q(X)_{\hat{x}})$ for every closed point $x\in X$, while each $\omega_x$ defines left Haar measure on $G(F_q(X)_{\hat{x}})$. Tamagawa measure $\omega_G$ on its adelic points $G({\mb A}_X)$ is defined as $q^{(1-g_X)d}\prod_{x\in X}L_x(\mc V_1,1)\omega_x$ and does not depend on the choice of $\omega$. Tamagawa number $\tau_G$ is defined as $\rho_{\widehat G}^{-1} \omega_G(G(K)\backslash G^1(\mb A))$. Quasi-discriminant $D_G$ denotes $\omega_G(G^c(\mb A_X))^{-2}$, where $G^c(\mb A_X)\subset G({\mb A}_X)$ is the maximal compact subgroup. 

It is easy to see that for $G$ split and $\omega$ defined over $\mb F_q$ we have 
$$
\omega_x(G(\mc O_x))=\frac{|G(\mb F_{q^{\deg x}})|}{q^{d\deg x}}=\prod_{i\geq1} \det(1-q^{-i\deg x}\mr{Fr}_x|V_i)= \prod_{i\geq1}L_x(\mc V_i, i)= \prod_{i\geq1}L_x(\mc V_i(1-i), 1)
$$
, so, as $G^c(\mb A_X)=\prod_{x\in X}G(\mc O_x)$ we get that 
$$
\omega_G(G^c(\mb A_X))=\prod_{i\geq 2} L(\mc V_i(1-i), 1)^{-1}= \prod_{i\geq 2} \rho_{\mc V_i(1-i)}^{-1}
$$.

The Siegel's mass formula (see \cite{HN})says that 
$$
 M_G(X)= \rho_{\widehat G}\tau_G D_G^{1/2}= \tau_Gq^{g_X-1}\dim G\prod_{i\geq 1} \rho_{\mc V_i(1-i)}
$$

Let's now take some asymptotically exact system $(\Xi, \ul{\mb Q_l})$. What is the asymptotic of $M_G(X_i)$ as $i$ grows? Applying the result of the previous section, we can easily give the answer.

\begin{thm}
Let $(\Xi, \ul{\mb Q_l})$ be asymptotically exact and let $G$ be a split reductive group over $\mb F_q$. Then

$$
\lim_{i\rightarrow\infty}\frac{\log_q\!M_G(X_i)}{g_{X_i}}= \dim G - \sum_{r=1}^{\infty}\gamma_{r}(\Xi, \ul{\mb Q_l})\log_q\!\left(\frac{|G(\mb F_{q^r})|}{q^{r\dim G}}\right)
$$
\begin{rem}
for $G=\mb G_m$, $\mr{Bun}^0_G = \mr{Pic}^0(X)$ and we get classical Tsfasman-Vl\v{a}duts formula:

$$
\lim_{i\rightarrow\infty}\frac{\log_q\!|\mr{Pic}^0(X_i)|}{g_{X_i}}= 1 - \sum_{r=1}^{\infty}\gamma_{r}(\Xi, \ul{\mb Q_l})\log_q\!\left(\frac{q^r -1}{q^r}\right)
$$

\end{rem}
\end{thm}

\textit{Proof:}

Applying $\log_q$ to Siegel formula we get

\begin{eqnarray*}
\lim_{i\rightarrow\infty}\frac{\log_q\!M_G(X_i)}{g_{X_i}}&=& \lim_{i\rightarrow\infty}\frac{\log_q\!\tau_G}{g_{X_i}} + \lim_{i\rightarrow\infty}\frac{(g_{X_i}-1)\dim G }{g_{X_i}}+  \lim_{i\rightarrow\infty}\frac{\log_q\!\prod_{i\geq 1} \rho_{\mc V_i(i-1)}}{g_{X_i}} =\\
&&=  \lim_{i\rightarrow\infty}\frac{\log_q\!\tau_G}{g_{X_i}} + \dim G -  \sum_{r=1}^{\infty}\gamma_{r}(\Xi, \ul{\mb Q_l})\log_q\!\left(\frac{|G(\mb F_{q^r})|}{q^{r\dim G}}\right),
\end{eqnarray*}

where the last equality is due to Steinberg's formula. Now, assuming Weil conjectures for semisimple groups over function fields (claimed to be proved by Lurie and Gaitsgory, see \cite{GL} for the first written down part of the proof) we can use [\cite{tamadawa}, Corollary 6.8.], which says that provided $G$ split, $\tau_G$ is constant (does not depend on $X_i$). So $\lim_{i\rightarrow\infty}\frac{\log_q\!\tau_G}{g_{X_i}} =0$ and we get the desired formula.

\subsection{Asymptotic formula for stable bundles}
 Here we restrict ourselves to the case $G=\mr{GL}_n$. We are going to extend the formula for $M_{\mr{GL}_n}(X)$ to its semistable part $M^{ss}_{\mr{GL}_n}(X)$. Now we do not have anything like Siegel mass formula, but it comes out that $M^{ss}_{\mr{GL}_n}(X)$ can be expressed through $M_{\mr{GL}_m}(X)$ for $m\le n$ (this is the formula of Don Zagier, see [\cite{M}, Theorem 5.3] or [\cite{Z}, Theorem 3]):
 
$$
M^{ss}_{\mr{GL}_n}(X)= \sum_{n_1+\cdots +n_k=n}q^{(g_X-1)\sum_{i<j}n_in_j}
\prod_{i=1}^{k-1}\frac{q^{(n_i+n_{i+1})\{(n_1+\cdots+n_i)n/d\}}}{1-q^{n_i+n_{i+1}}}M_{\mr{GL}_{n_1}}(X)\cdots M_{\mr{GL}_{n_k}}(X)
$$

We need to show that the first summand with $n_1=n, n_i =0$ is asymptotically the biggest. 

\begin{eqnarray}
\lim_{i\rightarrow\infty} \frac{\log_q M_{\mr{GL}_n}(X_i)}
{g_{X_i}} &=& n^2 - \sum_{r=1}^{\infty}\gamma_{r}(\Xi, \ul{\mb Q_l})\log_q\!\left(\frac{\prod_{i=1}^n (q^{rn}-q^{r(n-i)})}{q^{r n^2}}\right)=\\ &=& n^2 - \sum_{r=1}^{\infty}\gamma_{r}(\Xi, \ul{\mb Q_l})\sum_{i=1}^n \log_q\!\left(1-q^{-ir}\right).
\end{eqnarray}

Now, bounding $-\log_q\!\left(1-q^{-ir}\right)$ by geometric progression, we get that $\log_q\!\left(1-q^{-ir}\right)\le  q^{-ir+1}$, so again using Tsfasman-Vl\v{a}duts bound \cite{TV}

$$
\lim_{i\rightarrow\infty} \frac{\log_q M_{\mr{GL}_n}(X_i)}
{g_{X_i}}\le n^2 - \sum_{r=1}^{\infty}\gamma_{r}(\Xi, \ul{\mb Q_l}) n q^{-ir+1} < n^2 - n.
$$

Now, let's look at the limit of each summand: for $n_1,n_2,\ldots,n_k$ we get 
$$
\sum_{i<j}n_in_j + \sum_{j=1}^k \lim_{i\rightarrow\infty} \frac{\log_q M_{\mr{GL}_{j}(X_i)}}
{g_{X_i}}< \sum_{i<j}n_in_j + \sum_{j=1}^k n_j^2+n_j
$$

But as $n=\sum n_i$, $n^2=(\sum n_i)^2\ge \sum_{i<j}n_in_j + \sum_{j=1}^k n_j^2+n_j$. In logarithm term for $n_1=n$ we have $n^2+\mbox{something positive}$, because $|\mr{GL_n}(\mb F_q)|<q^{n^2}$, so it grows faster than any other, so 

$$
\lim_{i\rightarrow\infty} \frac{\log_q M_{\mr{GL}_n}(X_i)}
{g_{X_i}}=\lim_{i\rightarrow\infty} \frac{\log_q M^{ss}_{\mr{GL}_n}(X_i)}
{g_{X_i}}
$$

Let $\mc M^0_{\mr{GL}_n(X)}$ denote the moduli space of semistable vector bundles on $X$, whose determinant bundle is trivial. Now, as there are only finitely many semistable bundles and the order of automorphism group of semistable bundle is bounded by the order of $|\mr{GL}_n(\mb F_q)|$, the asymptotic of $M^{ss}_{\mr{GL}_n}(X_i)$ is the same as for $|\mc M^0_{\mr{GL}_n(X)}(\mb F_q)|$. So we get the following theorem

\begin{thm}
Let $(\Xi, \ul{\mb Q_l})$ be asymptotically exact. Then

$$
\lim_{i\rightarrow\infty}\frac{\log_q\!|\mc M^{ss}_{\mr{GL}_n(X)}(\mb F_q)|}{g_{X_i}}= \dim G - \sum_{r=1}^{\infty}\gamma_{r}(\Xi, \ul{\mb Q_l})\log_q\!\left(\frac{|G(\mb F_{q^r})|}{q^{r\dim G}}\right)
$$

\end{thm}

\medskip

\noindent {\bf Acknowledgements.} This is essentially the result of my first-year study as master at HSE and is almost identical to my coursework. I would like to thank my teachers Michael Finkelberg and Alexey Zykin.


\begin{thebibliography}{XXX}
\bibitem[BD]{tamadawa}K.~Behrend, A.~Dhillon, {\em Connected components of moduli stacks of torsors via Tamagawa numbers}, \arXiv{0503383v2}

\bibitem[HN]{HN} G.~Harder, M.~S.~Narasimhan, {\em On the cohomology groups of moduli spaces of vector
bundles on curves}, Math. Ann. 212 (1974/75), 215–248.

\bibitem[GL]{GL} D.~Gaitsgory, J.~Lurie, {\em Weil's conjecture for function fields I}, http://www.math.harvard.edu/~lurie/papers/tamagawa.pdf

\bibitem[M]{M} S.~Mozgovoy, {\em Poincare polynomials of moduli spaces of stable bundles over curves}, \arXiv{0711.0634v2}

\bibitem[TV]{TV} M.~A.~Tsfasman, S.~G.~Vl\v{a}duts. {\em Asymptotic properties of zeta-functions}, J.
Math. Sciences (New York), 1997, v.84, n.5, pp.1445-1467.

\bibitem[Z]{Z} D.~Zagier, {\em Elementary aspects of the Verlinde formula and of the Harder-Narasimhan-Atiyah-Bott formula}, Proceedings of the Hirzebruch 65 Conference on Algebraic Geometry
(Ramat Gan, 1993), Israel Math. Conf. Proc., vol. 9, Bar-Ilan Univ., 1996,
pp. 445–462.




\end{thebibliography}
\end{document}